\documentclass[11pt]{amsart}
\usepackage[margin=1.08in]{geometry}
\usepackage{amsmath,amssymb,amsthm,mathtools}
\usepackage[T1]{fontenc}
\usepackage{lmodern,microtype}
\usepackage[colorlinks=true,linkcolor=blue,citecolor=blue,urlcolor=blue]{hyperref}
\hypersetup{pdftitle={Unbounded total mean curvature with prescribed boundary metric}}
\newtheorem{theorem}{Theorem}[section]
\newtheorem{proposition}[theorem]{Proposition}
\newtheorem{conjecture}[theorem]{Conjecture}
\newtheorem{lemma}[theorem]{Lemma}

\theoremstyle{remark}

\newcommand{\Ric}{\operatorname{Ric}}

\newcommand{\tr}{\operatorname{tr}}
\newcommand{\Div}{\operatorname{div}}
\newcommand{\Hess}{\operatorname{Hess}}
\newcommand{\Vol}{\operatorname{Vol}}
\newcommand{\F}{\mathcal F}
\newcommand{\Lie}{\mathcal L}
\newcommand{\T}{\mathbb T}
\numberwithin{equation}{section}
\title[Unbounded total mean curvature]{Unbounded total mean curvature\\ with prescribed boundary metric}
\author{Zhixin Wang}
\address{Department of Mathematics, Shanghai Jiao Tong University, Shanghai, 201100}
\email{jhin@sjtu.edu.cn}
\date{October 26, 2026}
\begin{document}

\begin{abstract}
B\"ar has shown that a fixed boundary metric and a scalar-curvature
lower bound do not generally bound total mean curvature from above.
We extend this unboundedness to every smooth compact connected
manifold of dimension $n\ge3$ with nonempty boundary and every
prescribed boundary metric, using scalar-flat metrics.
While B\"ar's examples use warped products, our construction
uses a nonzero cubic variation of total mean curvature near
a flat product, followed by localization and conformal correction.
We also construct arbitrarily $C^k$-small perturbations of the
flat torical half-cylinder with nonnegative scalar curvature,
unchanged boundary metric, an exact product end, and positive
total boundary mean curvature.
\end{abstract}
\maketitle

\section{Introduction}

The total mean curvature of a boundary is closely related to
scalar curvature in the interior. A fundamental example is the
Brown--York inequality proved by Shi--Tam \cite{ST}. Let
$(\Omega^n,g)$ be a compact connected orientable Riemannian
manifold with nonnegative scalar curvature and positive outward
mean curvature on its boundary. If each boundary component
$(\Sigma_i,\gamma_i)$ admits an isometric embedding into
$\mathbb R^n$ as a strictly convex hypersurface, then
\[
 \int_{\Sigma_i}H_g\,d\sigma_{\gamma_i}
 \le
 \int_{\Sigma_i}H_0\,d\sigma_{\gamma_i},
\]
where $H_0$ is the mean curvature of the Euclidean embedding.
The formulation in arbitrary dimensions without a spin assumption
uses the positive mass theorem of Brendle--Wang \cite{BW}.
This inequality distinguishes two geometric conditions:
convexity of the reference embedding and positivity of the
mean curvature of the physical boundary.

More generally, a fill-in of a closed Riemannian manifold
$(\Sigma,\gamma)$ is a compact Riemannian manifold $(\Omega,g)$
whose boundary, with its induced metric, is isometric to
$(\Sigma,\gamma)$. We use the outward unit normal $\nu$ and
the convention
\[
 A(U,V)=g(\nabla_U\nu,V),\qquad H=\tr A.
\]
Thus the Euclidean unit sphere has mean curvature $n-1$.
We write
\begin{equation}\label{notation-F}
 \F(g)=\int_{\partial\Omega}H_g\,d\sigma_g.
\end{equation}
Gromov proposed an upper bound for total mean curvature
depending only on the boundary metric and a scalar-curvature
lower bound \cite{Gro}. The following formulation is stated
explicitly in \cite[Conjecture 1.1]{FHH} and
\cite[Introduction]{Bar}.

\begin{conjecture}
For every closed Riemannian manifold $(\Sigma,\gamma)$ and
$\sigma\in\mathbb R$, there is a finite constant
$C(\Sigma,\gamma,\sigma)$ such that every smooth fill-in
$(\Omega,g)$ with $R_g\ge\sigma$ satisfies
\[
 \int_\Sigma H_g\,d\sigma_\gamma
 \le C(\Sigma,\gamma,\sigma).
\]
\end{conjecture}

Without a sign condition or a common lower bound on $H_g$,
this formulation is false. B\"ar
\cite{BarLower} constructs scalar-flat
warped metrics on a fixed cylinder with fixed boundary metrics
and total mean curvature tending to infinity. He also obtains
examples with connected boundary and uniformly bounded scalar
curvature by joining the two boundary components with a handle. Frenck--Hanke--Hirsch \cite{FHH} also construct
examples with unbounded total mean curvature while controlling
the boundary volume and diameter, but the boundary metric
varies along their sequence.
We show that this unboundedness holds on every smooth compact
connected manifold of dimension $n\ge3$ with nonempty boundary,
for every prescribed boundary metric, while retaining
scalar-flatness.

\begin{theorem}\label{thm:universal}
Let $\Omega$ be any smooth compact connected $n$-manifold with
nonempty boundary, where $n\ge3$, and let $\gamma$ be any smooth
Riemannian metric on $\partial\Omega$. There are smooth metrics
$g_N$ on $\Omega$ such that
\[
 R(g_N)=0,\qquad
 g_N|_{T\partial\Omega}=\gamma,\qquad
 \F(g_N)\longrightarrow+\infty.
\]
\end{theorem}

In particular, the theorem applies to the ball with its standard
round boundary metric and to every compact connected filling
of a prescribed flat torus. 

The local mechanism of our construction differs from that
of B\"ar's scalar-flat warped-product examples.
In his scalar-flat cylinder examples, the two boundary metrics
are $\gamma$ and $2^{4/n}\gamma$: they are fixed along the
sequence but are not isometric. Our construction starts with
perturbations of a flat product that preserve the same metric
on both boundary components of a finite torical band.
The first two variations of total mean curvature vanish,
while the third variation is nonzero.

This behavior cannot occur within the corresponding
warped-product ansatz. Indeed, if $\gamma$ is flat and
\[
 g=dt^2+\varphi(t)^2\gamma,\qquad w=\varphi^{n/2},
\]
then
\[
 R_g=-\frac{4(n-1)}{n}\frac{w''}{w}.
\]
Thus $R_g\ge0$ implies $w''\le0$. If
$\varphi(0)=\varphi(1)=1$, the total mean curvature is
\[
 \F(g)=\frac{2(n-1)}{n}\Vol_\gamma(\T^{n-1})
 \bigl(w'(1)-w'(0)\bigr)\le0.
\]
Hence this ansatz cannot produce positive total mean curvature
while preserving the same flat metric on both boundary components.

The local mechanism also yields perturbations of a flat
half-cylinder that are arbitrarily small in any prescribed
finite $C^k$ norm. For this construction, let
\begin{equation}\label{background-metric}
 M=I\times\T^{n-1},\qquad g_0=dt^2+\gamma,
\end{equation}
where
\[
 \T^{n-1}=(\mathbb R/2\pi\mathbb Z)^{n-1},
 \qquad
 \gamma=dx_1^2+\cdots+dx_{n-1}^2,
\]
and $I=[0,1]$ or $[0,\infty)$, as specified below.
All $C^k$ norms on these products are taken with respect
to $g_0$, uniformly on $M$.

\begin{theorem}\label{thm:half}
Let $I=[0,\infty)$ in \eqref{background-metric}, with $n\ge3$.
For every integer $k\ge2$ and every $\varepsilon>0$, there
is a smooth metric $g$ on $M$ such that
\[
 R_g\ge0,\qquad
 g|_{T\partial M}=\gamma,\qquad
 \int_{\partial M}H_g\,d\sigma_\gamma>0,\qquad
 \|g-g_0\|_{C^k(g_0)}<\varepsilon.
\]
Moreover, $g=g_0$ for all sufficiently large $t$, and $H_g$
changes sign on $\partial M$. The support of the perturbation
is allowed to depend on $g$.
\end{theorem}

Theorem~\ref{thm:half} gives arbitrarily small perturbations
in any prescribed finite $C^k$ norm, with support allowed
to depend on the perturbation.
Theorem~\ref{thm:universal} instead concerns unbounded total
mean curvature on a fixed compact manifold.

The distinction between pointwise and integral control of
mean curvature is important here. Shi--Wang--Wei \cite{SWW}
prove that every smooth boundary metric on a compact manifold
extends to a metric of positive scalar curvature. They also
establish an obstruction to nonnegative scalar curvature
fill-ins when the prescribed mean curvature is sufficiently
large at every point, assuming that the boundary admits a
topological embedding into Euclidean space. Miao
\cite{MiaoLarge} removes this embedding hypothesis for
boundaries that admit a compact filling. These results bound
the minimum of $H$ from above; they do not give an upper bound for its integral.

A common lower bound on $H$ can nevertheless restore an
upper bound for its integral. Frenck--Hanke--Hirsch
\cite{FHH} and B\"ar \cite{Bar} establish such estimates for
spin fill-ins, allowing sign-changing mean curvature.
Consequently, on a spin filling, any sequence as in
Theorem~\ref{thm:universal} must satisfy
\[
 \min_{\partial\Omega}H_{g_N}\longrightarrow-\infty.
\]
In the half-cylinder setting, Lemma~\ref{lem:half-rigidity}
shows that nonnegative boundary mean curvature forces
the metric to split as a flat product. Thus the positive
total mean curvature in Theorem~\ref{thm:half} necessarily
comes with a change of sign.

Brown--York type inequalities also extend to non-Euclidean
reference geometries. Lu--Miao \cite{LuMiao} establish a
weighted Brown--York type inequality in the spatial
Schwarzschild manifold, incorporating the area of an
outermost minimal inner boundary through a localized
Penrose inequality. For toroidal reference geometries, the author
\cite{ZW} first establishes a positive mass type theorem for a
singular static model and then derives a Brown--York type
inequality for flat torus boundaries
\cite[Theorem 1.4]{ZW}. Alaee--Hung--Khuri \cite{AHK} establish
a weighted analogue in the Kottler setting. Together with
the Shi--Tam theorem \cite{ST}, these results provide sharp
comparisons of total, or suitably weighted total, mean
curvature with reference geometries, with additional
terms accounting for inner boundary or singular-end data
where appropriate.

The importance of strict reference convexity is illustrated
by the counterexamples of Brendle--Marques--Neves to Min-Oo's
conjecture \cite{BMN}. For the hemisphere, whose boundary is
totally geodesic, they construct metrics with the same induced
boundary metric, larger scalar curvature, and strictly
positive boundary mean curvature. In contrast,
Brendle--Marques \cite{BM} obtain local rigidity for
sufficiently small spherical caps through a coercive quadratic
estimate under scalar and boundary mean curvature comparison.
More recently, Zhang--Zhao \cite{ZZ} establish rigidity for
sufficiently small spherical caps among arbitrary oriented
fill-ins in every dimension, assuming the same boundary
metric and
\[
 R_g\ge n(n-1),\qquad H_g\ge(n-1)\cot r,
\]
without a closeness assumption.

These results illustrate two different roles of the
boundary hypotheses in Brown--York type inequalities.
Convexity of the reference boundary enters comparison
and local rigidity results, while a lower bound on
the physical boundary mean curvature may be needed
even for a finite upper bound on its integral.
The latter issue persists for a single fixed flat
torus boundary with nonnegative scalar curvature,
as our examples show.

Throughout the paper, an integral over $\partial M$ includes
both boundary components in the finite-band case and only
$\{0\}\times\T^{n-1}$ in the half-cylinder case.
We use $\Delta=\Div\nabla$ and write
\[
 q=\frac4{n-2},\qquad c_n=\frac{4(n-1)}{n-2}.
\]

\section{The cubic mechanism}
We outline the cubic mechanism on a finite flat band
and explain how it leads to the two main theorems.
The detailed half-cylinder construction is carried out
in the next section.

Our construction is motivated by the integral identity of
Brendle--Marques \cite{BM}. Let $I=[0,1]$ in \eqref{background-metric}. On our flat background, for
$g=g_0+sh$ with fixed induced boundary metric and
$\Div_{g_0}h=0$, the corresponding expansion reads
\[
 \int_M R_g\,dV_{g_0}
 +\int_{\partial M}(2-sh_{\nu\nu})H_g\,d\sigma_\gamma
 =
 -\frac{s^2}{4}\int_M
 \bigl(|\nabla h|^2+|\nabla\tr_{g_0}h|^2\bigr)\,dV_{g_0}
 +O(s^3).
\]
Here $\nu$ is the outward unit normal of $g_0$.
The above expansion is the flat-background specialization of
the weighted integral identity of Brendle--Marques \cite{BM},
which they use to prove local scalar-curvature rigidity for
sufficiently small geodesic balls in $S^n$. In our setting, however, this estimate alone does not yield rigidity. Indeed, consider a flat
annulus $(A,g_A)$ whose two boundary curves have the same
prescribed length but need not be geodesics. The product
\begin{equation}\label{motivation-metric}
 (A,g_A)\times(\T^{n-2},\gamma')
\end{equation}
has scalar curvature zero and the prescribed flat product
metric on both boundary components. By Gauss--Bonnet, its total boundary mean curvature is zero.

This example motivates examining infinitesimal boundary
bending. The variations below need not satisfy the
divergence-free condition in the preceding expansion;
we compute their boundary variations directly.
Choose $f\in C^\infty(\T^{n-1})$ and set
\begin{equation}\label{band-variation-def}
 X=tf\,\partial_t,\qquad
 g_s=u_s^{4/(n-2)}(g_0+s\Lie_Xg_0),
\end{equation}
where $u_s$ is chosen so that $R_{g_s}=0$ and
$u_s|_{\partial M}=1$. The induced metric on both boundary
components is unchanged, and $u_s=1+O(s^2)$. The initial tangent of the corrected family is
$h=\Lie_Xg_0$. By naturality of the Ricci tensor and flatness
of $g_0$, the first and second variations vanish, and the expansion begins at cubic order:
\[
 \int_{\partial M}H_{g_s}\,d\sigma_\gamma
 =
 -\frac{s^3}{6}\int_{\T^{n-1}}
 f\bigl((\Delta_\gamma f)^2-|\Hess_\gamma f|^2\bigr)
 \,d\sigma_\gamma
 +O(s^4).
\]

The boundary bending suggested by the product example
\eqref{motivation-metric} depends on only one torus direction and the cubic term vanishes. This agrees with the Gauss--Bonnet
identity. To obtain a nonzero cubic
term, we instead choose
\[
 f(x)=\cos x_1+\cos x_2+\cos(x_1+x_2).
\]
Thus, for sufficiently small negative $s$, the metric $g_s$ on
$[0,1]\times\T^{n-1}$ is scalar-flat, induces the same flat
metric $\gamma$ on both boundary components, and satisfies
\[
 \int_{\partial M}H_{g_s}\,d\sigma_\gamma>0.
\]
This nonzero cubic variation is the central mechanism
underlying our construction.

To obtain a sequence of metrics with arbitrarily large total
mean curvature, we adapt this mechanism to the half-cylinder
to prove Theorem~\ref{thm:half}, and let $g$ be a metric
provided by that theorem.
To amplify its total mean curvature, for each positive
integer $N$ define
\[
 \Phi_N(t,x)=(Nt,Nx),\qquad
 g_N=N^{-2}\Phi_N^*g,
\]
where multiplication by $N$ is understood modulo $2\pi$
in each torus coordinate. Then
\[
 R_{g_N}=N^2R_g\circ\Phi_N\ge0,\qquad
 g_N|_{T\partial M}=\gamma,\qquad
 \F(g_N)=N\F(g).
\]

To prove Theorem~\ref{thm:universal}, we localize the same
high-frequency construction in a boundary coordinate chart.
A smooth family of periodic half-cylinder metrics allows
the model to vary with the prescribed boundary metric
and to return to the background metric outside the chart.
The resulting metrics preserve the boundary metric exactly
and have total mean curvature at least $aN-C$, where $a>0$.
Although localization may introduce negative scalar
curvature, its magnitude is at most $CN$ and its support
lies in a boundary layer of thickness $O(N^{-1})$.
Lemma~\ref{lem:repair} then provides a scalar-flat conformal
correction preserving the boundary metric and decreasing
total mean curvature by at most a constant independent of
$N$. The corrected metrics therefore satisfy
$\F\to+\infty$.

\section{Local deformation of the half-cylinder }
In this section, consider $I=\mathbb{R}_+$ in \eqref{background-metric}. The proof proceeds in three steps. First, as in the previous
section, we start with a small, compactly supported perturbation
in a Lie derivative direction and apply a conformal correction
to obtain a scalar-flat metric with the prescribed boundary
metric. Here we also need estimates for the conformal factor,
including its limiting value and decay at infinity. Second,
we use these estimates to extend the third-variation calculation
to the half-cylinder and obtain positive total mean curvature.
Finally, a cutoff construction followed by a further conformal
change restores the exact product metric $g_0$ near infinity,
while preserving nonnegative scalar curvature and positivity
of the boundary integral. This proves Theorem~\ref{thm:half}.

Choose $\chi\in C^\infty([0,\infty))$ equal to $1$ near zero and
$0$ for $t\ge1$. Let $f\in C^\infty(\T^{n-1})$ and regard it as independent of $t$. Set
\begin{equation}\label{cylinder-variation-set}
    \psi=\chi f,\quad X=-\psi\partial_t,\quad \widetilde g_s=g_0+s\Lie_Xg_0
\end{equation}
In product coordinates,
\begin{equation}\label{eq:seed}
 (\widetilde g_s)_{tt}=1-2s\chi'f,\qquad
 (\widetilde g_s)_{ta}=-s\chi f_a,\qquad
 (\widetilde g_s)_{ab}=\gamma_{ab}.
\end{equation}
This is positive definite for small $s$, induces exactly $\gamma$
on every slice, and equals $g_0$ for $t\ge1$.

\subsection{Positive total mean curvature with nonnegative scalar curvature}
We need the formula for scalar curvature and mean curvature under the conformal change $ g_u=u^{4/(n-2)}g$
\begin{equation}\label{eq-conformal-RH}
\begin{aligned}
 R_{g_u}
 &=u^{-\frac{n+2}{n-2}}
 \left(-\frac{4(n-1)}{n-2}\Delta_g u+R_g u\right)\coloneqq u^{-\frac{n+2}{n-2}} L_g(u),\\
 H_{g_u}
 &=u^{-\frac{n}{n-2}}
 \left(\frac{2(n-1)}{n-2}\partial_{\nu_g}u+H_g u\right).
\end{aligned}
\end{equation}
\begin{lemma}\label{lem:solve}
For small $|s|$, the problem
\begin{equation}\label{eq-infiniteconf}
 -c_n\Delta_{\widetilde g_s}u_s+R(\widetilde g_s)u_s=0,
 \qquad u_s|_{t=0}=1,\qquad u_s\text{ bounded}
\end{equation}
has a unique bounded smooth solution. It is positive and smooth in $s$.
There exists a smooth function $c(s)>0$ and $\lambda>0$ such that, for
all fixed $j$,
\begin{equation}\label{eq-decay}
 \|u_s-1\|_{C^j(M)}\le C_js^2,\qquad
 |\nabla^j(u_s-c(s))|\le C_js^2e^{-\lambda t}\quad(t\ge2).
\end{equation}
Thus $g_s=u_s^q\widetilde g_s$ is scalar-flat, induces $\gamma$ at
the compact boundary, and approaches $c(s)^qg_0$ exponentially.
\end{lemma}
\begin{proof}
Write $L_s$ and $L_0$ for $L_{\widetilde g_s}$ and $L_{g_0}$ respectively. We first establish a bounded inverse for $L_0$ on sources
supported in the fixed cylinder
$K=[0,2]\times\T^{n-1}$.

For $F\in C^{0,\alpha}(M)$ supported in $K$, solve
\[
 L_0w_T=F
 \quad\text{on }[0,T]\times\T^{n-1},
 \qquad w_T(0)=w_T(T)=0,
 \qquad T>3.
\]
Existence and regularity follow from standard Dirichlet
theory; see \cite[Chapter 6]{GT}.
Choose a smooth nonnegative function $\rho(t)$, equal to
$1$ on $[0,2]$ and vanishing for $t\ge3$, and put
\[
 b(t)=\frac1{c_n}\int_0^\infty\min\{t,r\}\rho(r)\,dr.
\]
Then $b$ is bounded, nonnegative, and
$L_0b=\rho$. The maximum principle gives $|w_T|\le \|F\|_{C^0}b$.
Letting $T\rightarrow \infty$, local and boundary Schauder estimates
therefore yield a bounded solution $w$ satisfying
\[
 L_0w=F,\qquad w|_{t=0}=0,\qquad
 \|w\|_{C^{2,\alpha}(M)}
 \le C\|F\|_{C^{0,\alpha}(M)}.
\]
The constant is independent of $F$. Uniqueness follows again
from the maximum principle: if $L_0w=0$, $w|_{t=0}=0$, and
$|w|\le B$, comparison with $\pm Bt/T$ on $[0,T]\times\T^{n-1}$
and passage to $T\to\infty$ give $w=0$. Denote this operator by $w=G(F)$.

Now set $P_s=L_s-L_0$. Its coefficients are supported in $K$,
and smooth dependence of $\widetilde g_s$ gives
\[
 \|P_sw\|_{C^{0,\alpha}(M)}
 \le C|s|\|w\|_{C^{2,\alpha}(M)}.
\]
Writing $u_s=1+w_s$, we obtain from \eqref{eq-infiniteconf}
\[
 (I+GP_s)w_s=-G R(\widetilde g_s).
\]
For small $|s|$, the operator $I+GP_s$ is invertible. This gives a unique bounded solution,
depending smoothly on $s$. Naturality implies
\[
 \left.\frac{d}{ds}\right|_{s=0}R(\widetilde g_s)
 =DR_{g_0}(\Lie_Xg_0)=0.
\]
Consequently, $|R(\widetilde g_s)|\leq Cs^2$
\[
 \|u_s-1\|_{C^{2,\alpha}(M)}\le Cs^2,
\]
and elliptic regularity gives
$\|u_s-1\|_{C^j(M)}\le C_js^2$ for every fixed $j$.
In particular, $u_s>0$ for small $|s|$.
Any bounded smooth solution satisfies the same global
regularity estimates and hence the same inverse equation,
so uniqueness holds among all bounded smooth solutions.

On $t\ge2$, the equation reduces to
$ (\partial_t^2+\Delta_\gamma)u_s=0$.
Expanding in eigenfunctions
$-\Delta_\gamma\phi_\ell=\lambda_\ell\phi_\ell$,
boundedness eliminates the term linear in $t$ associated
with the zero eigenvalue, as well as all exponentially
growing terms associated with positive eigenvalues. Thus
\[
 u_s(t,x)
 =c(s)+\sum_{\ell\ge1}b_\ell(s)
 e^{-\sqrt{\lambda_\ell}(t-2)}\phi_\ell(x).
\]
The preceding estimates imply \eqref{eq-decay} for any fixed
$0<\lambda<\sqrt{\lambda_1}$. Moreover,
\[
 c(s)=\frac{1}{\Vol_\gamma(\T^{n-1})}
 \int_{\T^{n-1}}u_s(2,x)\,d\sigma_\gamma,
\]
so $c(s)$ depends smoothly on $s$ and satisfies
$c(s)=1+O(s^2)>0$. The exponential decay also holds after
differentiation in $s$: for every fixed $j,\ell\ge0$,
\[
 \bigl|\nabla^j\partial_s^\ell(u_s-c(s))\bigr|
 \le C_{j,\ell}e^{-\lambda t}
 \qquad(t\ge2),
\]
uniformly for sufficiently small $|s|$.
The conformal transformation formula proves the remaining
assertions.
\end{proof}

\begin{lemma}\label{lem-cylinder}
The family of Lemma~\ref{lem:solve} satisfies
\begin{equation}\label{eq:cubic-expansion-cylinder}
 \int_{\partial M}H_{g_s}\,d\sigma_\gamma
 =
 -\frac{s^3}{6}\int_{\T^{n-1}}
 f\bigl((\Delta_\gamma f)^2-|\Hess_\gamma f|^2\bigr)
 \,d\sigma_\gamma
 +O(s^4),
\end{equation}
with a remainder uniform for small $s$.
\end{lemma}

\begin{proof}
Write $\Sigma=\{0\}\times\T^{n-1}$ and
\[
 \F(s)=\int_\Sigma H_{g_s}\,d\sigma_\gamma.
\]
By Lemma~\ref{lem:solve}, $u_0=1$ and $\partial_su_s|_{s=0}=0$, so
\[
 \partial_sg_s|_{s=0}=\Lie_Xg_0,\qquad
 X=-\chi f\partial_t,\qquad X|_\Sigma=f\nu,
\]
where $\nu=-\partial_t$ is the outward unit normal of $g_0$.

We first justify the first variation formula on the
half-cylinder. Set $v=\partial_sg_s$, with all geometric
quantities below computed with respect to $g_s$.
The scalar-curvature variation is
\begin{equation}\label{R-1st-variation-formula}
 R'_s
 =
 -\langle\Ric(g_s),v\rangle
 +\Div\Div v-\Delta\tr v.
\end{equation}
On a fixed boundary hypersurface, the mean-curvature
variation satisfies
\[
 2H'_s
 =
 \bigl(d(\tr v)-\Div v\bigr)(\nu_s)
 -\Div_{\partial M}W-\langle A_s,v^T\rangle,
 \qquad
 W=\bigl(v(\nu_s,\cdot)|_{T\partial M}\bigr)^\sharp.
\]
Here $v^T$ denotes the tangential restriction of $v$.

Apply these formulas on
$M_T=[0,T]\times\T^{n-1}$, with
$\Sigma_T=\{T\}\times\T^{n-1}$.
Since $R(g_s)=0$ and $v^T=0$ on $\Sigma$, integration
and the divergence theorem give
\[
 2\F'(s)
 =
 -\int_{M_T}\langle\Ric(g_s),v\rangle\,dV_{g_s}
 -\int_{\Sigma_T}
   \bigl(2H'_s+\langle A_s,v^T\rangle\bigr)
   \,d\sigma_{g_s}.
\]
On the product end, $g_s=u_s^qg_0$, and $u_s-c(s)$
decays exponentially, together with its spatial and
parameter derivatives.
Consequently, $A_s$, $H'_s$, and their parameter
derivatives decay exponentially on $\Sigma_T$, whereas
$v$ and the slice volumes remain uniformly bounded.
The integral over $\Sigma_T$ therefore tends to zero.
The same decay makes the interior integral convergent
and permits differentiation twice in $s$ under the
integral sign. Hence
\[
 2\F'(s)
 =
 -\int_M
 \langle\Ric(g_s),\partial_sg_s\rangle_{g_s}\,dV_{g_s}.
\]

Flatness and naturality give
\[
 \Ric(g_0)=0,\qquad
 \left.\partial_s\Ric(g_s)\right|_{s=0}
 =D\Ric_{g_0}(\Lie_Xg_0)
 =\Lie_X\Ric(g_0)=0.
\]
It follows that $\F'(0)=\F''(0)=0$.
Set $
 T=\left.\partial_s^2\Ric(g_s)\right|_{s=0}$.
Differentiating the first variation formula twice yields
\[
 \F'''(0)
 =
 -\frac12\int_M
 \langle T,\Lie_Xg_0\rangle_{g_0}\,dV_{g_0}.
\]
The twice differentiated contracted Bianchi identity,
together with $R(g_s)=0$ and the vanishing of
$\Ric(g_0)$ and its first variation, gives
$\Div_{g_0}T=0$.
Since $X$ is compactly supported, integration by parts
has no boundary contribution at infinity and gives
\begin{equation}\label{eq:third}
 \F'''(0)
 =
 -\int_\Sigma T(X,\nu)\,d\sigma_\gamma
 =
 -\int_\Sigma fT(\nu,\nu)\,d\sigma_\gamma.
\end{equation}

The induced metric on $\Sigma$ remains flat.
The Gauss equation and scalar flatness therefore give
\[
 2\Ric(g_s)(\nu_s,\nu_s)=H_s^2-|A_s|^2.
\]
At $s=0$, we have $A_0=H_0=0$, and both the Ricci
tensor and its first variation vanish.
Differentiating twice gives
\[
 T(\nu,\nu)=(H'_0)^2-|A'_0|^2.
\]
Since $g'_0=\Lie_Xg_0$ and $X|_\Sigma=f\nu$,
the normal variation formulas on the flat product yield
\[
 A'_0=-\Hess_\gamma f,\qquad
 H'_0=-\Delta_\gamma f.
\]
Substitution into \eqref{eq:third} gives
\[
 \F'''(0)
 =
 -\int_{\T^{n-1}}
 f\bigl((\Delta_\gamma f)^2-|\Hess_\gamma f|^2\bigr)
 \,d\sigma_\gamma.
\]
Finally, $\F(0)=0$, and $\F(s)$ is smooth because
$\Sigma$ is compact and $g_s$ depends smoothly on $s$.
Taylor's theorem proves
\eqref{eq:cubic-expansion-cylinder}, with a uniform
$O(s^4)$ remainder on a sufficiently small parameter
interval.
\end{proof}

\subsection{Restoring the product end}
The scalar-flat family $g_s$ approaches $c(s)^qg_0$
at infinity. We restore the original product end while
controlling the loss in total mean curvature, so that
the cubic gain can be retained once $f$ and $s$ are chosen.
We restore the original product end with a loss that can
be made smaller than the cubic gain.
To do this, we first show $c(s)<1$, allowing
the restoring conformal factor to increase from its
boundary value while remaining superharmonic.
\begin{lemma}\label{lem-c-estimate}
Let $E=\int_{\T^{n-1}}|df|^2\,d\sigma_\gamma$ and write $V=\Vol_\gamma(\T^{n-1})$. Then
\begin{equation}\label{eq:cscale}
 c(s)=1-\frac{E}{c_nV}s^2+O(s^3).
\end{equation}
In particular $c(s)<1$ for small nonzero $s$ if $f$ is nonconstant.
\end{lemma}
\begin{proof}
Let $\widehat M=\mathbb R\times\T^{n-1}$ carry the flat product
metric $\widehat g_0=d\tau^2+\gamma$, and define
\[
 \Phi_s:M\longrightarrow\widehat M,
 \qquad \Phi_s(t,x)=(t-s\psi(t,x),x).
\]
For small $|s|$, this map is a local diffeomorphism. Pulling
back the target metric $\widehat g_0$ gives
\[
 g_s^{\rm flat}:=\Phi_s^*\widehat g_0
 =(dt-s\,d\psi)^2+\gamma
 =\widetilde g_s+s^2d\psi\otimes d\psi.
\]
where the last equality follows from the definition \eqref{cylinder-variation-set} and $ \Lie_Xg_0=-dt\otimes d\psi-d\psi\otimes dt$. In particular, $g_s^{\rm flat}$ is flat.
Its boundary need not be fixed. Let $h=\mathcal{L}_Xg_0$. Since both families are scalar-flat and have initial derivative
$h$, the chain rule gives
\[
 \begin{aligned}
 0
 &=\left.\partial_s^2R(g_s)\right|_{s=0}=DR_{g_0}\!\left(\left.\partial_s^2g_s\right|_{s=0}\right)
   +D^2R_{g_0}(h,h),\\
 0
 &=\left.\partial_s^2R(g_s^{\rm flat})\right|_{s=0}=DR_{g_0}(2d\psi\otimes d\psi)
   +D^2R_{g_0}(h,h).
 \end{aligned}
\]
Subtracting these identities yields
\[
 DR_{g_0}\!\left(
 \left.\partial_s^2g_s\right|_{s=0}
 -2d\psi\otimes d\psi
 \right)=0.
\] Set $ Q:=\left.\partial_s^2g_s\right|_{s=0}
       -2d\psi\otimes d\psi$.
Then $DR_{g_0}(Q)=0$.
Integrating \eqref{R-1st-variation-formula} over
a torus slice eliminates the tangential divergence
terms, while the two terms involving
$\partial_t^2Q_{tt}$ cancel. Thus
\[
 0=-\frac{d^2}{dt^2}
 \int_{\T^{n-1}}\tr_\gamma Q^T\,d\sigma_\gamma.
\]
where $Q^T$ is the restriction to the torus direction.
The integral is therefore affine in $t$; boundedness makes
it constant.

At $t=0$, the fixed boundary metric gives
$\left.\partial_s^2g_s^T\right|_{s=0}=0$, and $\psi=f$.
Hence this constant is $-2E$, where
\[
 E=\int_{\T^{n-1}}|df|_\gamma^2\,d\sigma_\gamma.
\]
At infinity, $\psi=0$ and $g_s\to c(s)^qg_0$ with parameter
derivatives. Since $c(0)=1$ and
$\left.\partial_sc(s)\right|_{s=0}=0$, we obtain
\[
 Q\longrightarrow
 q\left.\partial_s^2c(s)\right|_{s=0}g_0.
\]
 Using
$(n-1)q=c_n$, comparison of the two limits gives
\[
 -2E=c_nV\left.\partial_s^2c(s)\right|_{s=0}.
\]
and the lemma follows.
\end{proof}

We next construct a metric $\hat g_{s,L}$ that agrees exactly with $g_0$ for sufficiently large $t$.
\begin{lemma}\label{lem:restore}
Fix a nonconstant $f$ and the family $g_s$ above. There exist $s_0>0$
and $L_0>3$ such that, for $0<|s|<s_0$ and $L\ge L_0$, there is a
smooth metric $\widehat g_{s,L}$ satisfying
\[
 R(\widehat g_{s,L})\ge0,\qquad
 \widehat g_{s,L}|_{T\partial M}=\gamma,\qquad
 \widehat g_{s,L}=g_0\quad(t\ge L+3).
\]
Moreover,
\begin{equation}\label{eq:cost}
 0\le\int_{\partial M}H_{g_s}\,d\sigma_\gamma-\int_{\partial M}H_{\hat{g}_{s,L}}\,d\sigma_\gamma
 \le C\left(\frac{s^2}{L}+s^2e^{-\lambda L}\right).
\end{equation}
 For each fixed $k$,
\begin{equation}\label{eq:Ck}
 \|\widehat g_{s,L}-g_0\|_{C^k(g_0)}\le C_k|s|.
\end{equation}
The constants are independent of $s,L$ in these ranges.

\end{lemma}
\begin{proof}
In this proof all the constants $C,C_j$ are independent of $s,L$. Cut off $g_s-c(s)^qg_0$ over $[L,L+1]$ to obtain a metric
$\bar g$ equal to $g_s$ for $t\le L$ and to $c(s)^qg_0$
for $t\ge L+1$. Let $\lambda$ be as in \eqref{eq-decay}, and we have
\[ \|\bar g-g_s\|_{C^j}
 +\|R(\bar g)\|_{C^j}
 \le C_js^2e^{-\lambda L},
 \qquad
 \operatorname{supp}R(\bar g)\subset[L,L+1].
\]
We seek a conformal correction $\widehat g=v(t)^q\bar g$,
where $v(0)=1$ and $v=c(s)^{-1}$ for $t\ge L+3$. Put
\begin{equation}\label{cylinder-B-alpha-def}
 \alpha=|dt|_{\bar g}^2,\qquad
 B=\frac{\Delta_{\bar g}t}{\alpha}.
\end{equation}
Since $B\leq Ce^{-\lambda t}$ for $2\le t\le L+1$ and vanishes for $t>L+1$, we can choose a smooth nonnegative function $a(t)\ge\sup_x B(t,x)$,
supported in $[0,L+2]$ with
\[
 \int_0^\infty a\,dt\le C,\qquad
 \|a\|_{C^j}\le C_j,
\]
and a nonnegative smooth bump $b_0$ supported in
$(L-1,L+2)$ satisfying 
\[
 c_n\alpha b_0\ge2R(\bar g)_-,
 \qquad
 \|b_0\|_{C^j}\le C_js^2e^{-\lambda L}.
\]
Let $\eta_L$ be a nonnegative unit-integral smooth bump
supported in $(L+2,L+3)$. For $\beta\ge0$, define
\begin{equation}\label{cylinder-y-def}
 y_\beta(t)=
 \int_t^\infty
 \exp\left(\int_t^r a(\rho)\,d\rho\right)
 \bigl(b_0(r)+\beta\eta_L(r)\bigr)\,dr.
\end{equation}
Then $y_\beta\ge0$ and
\begin{equation}\label{cylinder-y-ODE}
 y_\beta'+ay_\beta=-b_0-\beta\eta_L.
\end{equation}
We choose $\beta$ so that
\[
 \int_0^\infty y_\beta(t)\,dt=c(s)^{-1}-1.
\]
Indeed, $\int_0^{\infty}y_\beta(t)dt$ can be viewed as an affine function of $\beta$ with value
at most $CLs^2e^{-\lambda L}$ for $\beta=0$ and slope lying between $L+2$ and $C(L+3)$.
By Lemma~\ref{lem-c-estimate}, $\frac{1}{C}s^2<c(s)^{-1}-1<Cs^2$ and this determines a positive $\beta$ for
all sufficiently large $L$, uniformly for small nonzero $s$.
Moreover, $ \beta\le \frac{Cs^2}{L}$.
Set
\[
 v(t)=1+\int_0^t y_\beta(r)\,dr.
\]
Then
\[
 v(0)=1,\qquad 1\le v\le c(s)^{-1}<2,\qquad
 v=c(s)^{-1}\quad(t\ge L+3).
\]
Then we have
\begin{equation}
\begin{aligned}
     \Delta_{\bar g}v
 &=\alpha(y_\beta'+By_\beta)\\
 &\leq \alpha(y_\beta'+ay_\beta)\\
 &=\alpha(-b_0-\beta \eta_L)
 \le-\alpha b_0,
\end{aligned}
\end{equation}
where we used the chain rule and \eqref{cylinder-B-alpha-def} for the first line, $B\leq a$ for the second line and \eqref{cylinder-y-ODE} in the last line. And hence
\[
 -c_n\Delta_{\bar g}v+R(\bar g)v
 \ge c_n\alpha b_0-2R(\bar g)_-\ge0.
\]
Thus $\widehat g_{s,L}=v^q\bar g$ has nonnegative scalar
curvature by \eqref{eq-conformal-RH}, induces $\gamma$ at $\{0\}\times T^{n-1}$, and equals
$g_0$ for $t\ge L+3$. 

It remains to estimate the change in total mean curvature.
At $t=0$, we have $\bar g=g_s$ and $v=1$. The conformal
mean-curvature formula \eqref{eq-conformal-RH} gives
\[
 H_{\widehat g_{s,L}}-H_{g_s}
 =\frac{c_n}{2}\partial_{\nu_{g_s}}v
 =-\frac{c_n}{2}y_\beta(0)|dt|_{g_s},
\]
where $\nu_{g_s}$ is the outward unit normal.
By \eqref{cylinder-y-def}, the uniform bound on
$\int_0^\infty a\,dt$, and the estimates for $b_0$ and $\beta$,
\[
 0\le y_\beta(0)
 \le C\bigl(\|b_0\|_{L^1}+\beta\bigr)
 \le C\left(s^2e^{-\lambda L}+\frac{s^2}{L}\right).
\]
Since $|dt|_{g_s}$ is uniformly bounded on the fixed boundary,
integration yields
\[
 0\le
 \int_{\partial M}H_{g_s}\,d\sigma_\gamma
 -
 \int_{\partial M}H_{\widehat g_{s,L}}\,d\sigma_\gamma
 \le C\left(\frac{s^2}{L}+s^2e^{-\lambda L}\right),
\]
which proves \eqref{eq:cost}. Finally, the ODE for $y_\beta$ and the uniform derivative
bounds for $a,b_0,\eta_L$ give
\[
 \|v-1\|_{C^k}\le C_ks^2.
\]
Combining this with the cutoff estimate and
$\|g_s-g_0\|_{C^k}\le C_k|s|$ proves \eqref{eq:Ck}.

\end{proof}

\subsection{Proof of Theorem~\ref{thm:half}}

We now choose $f$ and combine the cubic expansion with
the product-end restoration.

\begin{proof}[Proof of Theorem~\ref{thm:half}]
Choose
\begin{equation}\label{choice-f}
 f(x)=\cos x_1+\cos x_2+\cos(x_1+x_2)
 \coloneqq a+b+c.
\end{equation}
A direct calculation gives
\[
 (\Delta_\gamma f)^2-|\Hess_\gamma f|^2
 =2(ab+ac+bc).
\]
The torus average of $abc$ is $1/4$, whereas
$a^2b,a^2c,ab^2,b^2c,ac^2,bc^2$ all have zero average.
Consequently,
\[
 \int_{\T^{n-1}}
 f\bigl((\Delta_\gamma f)^2-|\Hess_\gamma f|^2\bigr)
 \,d\sigma_\gamma
 =\frac32\Vol_\gamma(\T^{n-1}).
\]
Lemma~\ref{lem-cylinder} therefore gives
\[
 \int_{\partial M}H_{g_s}\,d\sigma_\gamma
 =-\frac14\Vol_\gamma(\T^{n-1})s^3+O(s^4).
\]
Take $s<0$ sufficiently small and set $L=|s|^{-2}$.
By Lemma~\ref{lem:restore}, the metric $\widehat g_{s,L}$
has nonnegative scalar curvature, induces $\gamma$ on
$\partial M$, and equals $g_0$ for $t\ge L+3$.
Since \eqref{eq:cost} gives a loss of order $O(s^4)$,
\[
 \int_{\partial M}H_{\widehat g_{s,L}}\,d\sigma_\gamma
 =\frac14\Vol_\gamma(\T^{n-1})|s|^3+O(s^4)>0
\]
for sufficiently small $|s|$.
For any fixed $k\ge2$ and $\varepsilon>0$, estimate
\eqref{eq:Ck} ensures that
\[
 \|\widehat g_{s,L}-g_0\|_{C^k(g_0)}<\varepsilon
\]
after decreasing $|s|$ further.

Since the boundary integral of $H_{\widehat g_{s,L}}$ is strictly positive, $H_{\widehat g_{s,L}}$ changes sign by Lemma~\ref{lem:half-rigidity} below. Taking $g=\widehat g_{s,L}$ completes the proof.
\end{proof}

\begin{lemma}\label{lem:half-rigidity}
Let $\gamma$ be a flat metric on $\T^{n-1}$, $n\ge3$,
and let $g$ be a smooth metric on
$M=[0,\infty)\times\T^{n-1}$ satisfying
\[
 R_g\ge0,\qquad
 g|_{T\partial M}=\gamma,\qquad
 g=dt^2+\gamma\quad\text{for all sufficiently large }t.
\]
If $H_g\ge0$ on $\partial M$, then $(M,g)$ is
isometric to the flat product
$([0,\infty)\times\T^{n-1},dt^2+\gamma)$.
In particular, $H_g\equiv0$.
\end{lemma}

\begin{proof}
Choose $T$ in the product end and truncate to
$M_T=[0,T]\times\T^{n-1}$.
Both boundary components have nonnegative outward
mean curvature. By the scalar-curvature rigidity
result of B\"ar--Hanke \cite[Theorem~2.19]{BH},
applied to the torical band $M_T$, we obtain
\[
 \Ric_g=0 \quad\text{on }M_T,\qquad
 H_g=0 \quad\text{on }\partial M_T.
\]

Let $u$ be harmonic on $M_T$, with boundary values
$0$ and $1$ on its two boundary components.
Reilly's formula, using $\Ric_g=0$, $H_g=0$,
and the constant boundary values, gives
\[
 \int_{M_T}|\Hess_g u|^2\,dV_g=0.
\]
Thus $\nabla u$ is a nonzero parallel vector field.
Its flow identifies $M_T$ isometrically with
$[0,\ell]\times(\T^{n-1},\gamma)$ for some $\ell>0$.
Attaching the remaining product end proves the claim.
\end{proof}

\section{Universal Result}
To localize the half-cylinder construction, we need models
that depend smoothly on the boundary metric and share
a common product end.

\begin{lemma}\label{lem-periodic-family}
Identify
$\T^{n-1}=\mathbb R^{n-1}/(2\pi\mathbb Z)^{n-1}$ and take
\[
 f(y)=\cos y_1+\cos y_2+\cos(y_1+y_2).
\]
For every positive definite matrix $Q_0$, there exist a
relatively compact neighborhood $\mathcal U$ of $Q_0$
in the positive definite cone, constants $s_*,A,b>0$,
and a smooth family of metrics $K_{Q,s}$ on $\mathbb{R}_+\times T^{n-1}$, for
$Q\in\mathcal U$ and $-s_*\le s\le0$, such that
\begin{align}
 &R(K_{Q,s})\ge0,\qquad
 K_{Q,s}|_{T\partial M}=Q,
 \label{eq:family-boundary}\\
 &K_{Q,0}=dt^2+Q,\qquad
 K_{Q,s}=dt^2+Q\quad(t\ge A),
 \label{eq:family-end}
\end{align}
and
\begin{equation}\label{eq:family-gain}
 \int_{\partial M}H_{K_{Q,-s_*}}\,d\sigma_Q\ge b
 \qquad(Q\in\mathcal U).
\end{equation}
All coefficient derivatives of fixed finite order,
including parameter derivatives, are uniformly bounded
on the indicated parameter set.
\end{lemma}
\begin{proof}
Carry out the construction with $\gamma=Q$ on a compact
matrix neighborhood of $Q_0$, and choose $\mathcal U$
with closure contained in its interior.

The construction in Lemma~\ref{lem:solve} depends smoothly
on the flat boundary metric $Q$.
Indeed, let $G_Q$ denote the bounded inverse of
$-c_n(\partial_t^2+\Delta_Q)$ on sources supported in
$[0,2]\times\T^{n-1}$, with zero boundary value at $t=0$.
The barrier used in that proof is independent of $Q$,
and the Schauder estimates are uniform on a compact
matrix neighborhood.

Smooth dependence of $G_Q$ can be seen directly in
the fixed Fourier basis. For $k\ne0$, put
$\mu_Q(k)=\sqrt{k^{\mathsf T}Q^{-1}k}$.
The Green kernel for the equation satisfied
by the $k$-th Fourier coefficient is
\[
 \frac{
 e^{-\mu_Q(k)|t-r|}
 -
 e^{-\mu_Q(k)(t+r)}
 }{2c_n\mu_Q(k)},
\]
For $k=0$, the Green kernel is $\min\{t,r\}/c_n$,
independent of $Q$.
Differentiation of these kernels, together with elliptic
estimates, gives smooth dependence on $Q$ and uniform
bounds for smooth sources supported in the fixed cylinder.

Writing
\[
 P_{Q,s}
 =L_{\widetilde g_{Q,s}}
   +c_n(\partial_t^2+\Delta_Q),
\]
the same inverse equation as in Lemma~\ref{lem:solve},
\[
 (I+G_QP_{Q,s})(u_{Q,s}-1)
 =-G_QR(\widetilde g_{Q,s}),
\]
therefore gives smooth dependence on $(Q,s)$ for uniformly
small $|s|$.
The uniform positive lower bound for the nonzero
eigenvalues of $-\Delta_Q$ also gives uniform exponential
decay, including parameter derivatives after decreasing
the decay rate if necessary.

Write $g_{Q,s}$ and $c(Q,s)$ for the resulting
scalar-flat family and limiting conformal factor.
Lemmas~\ref{lem-cylinder} and~\ref{lem-c-estimate} 
apply uniformly in $Q$.
The quadratic coefficient of $1-c(Q,s)$ has a uniform
positive lower bound. Moreover, the Fourier calculation
for the chosen $f$ gives
\[
 \int_{\T^{n-1}}
 f\bigl((\Delta_Qf)^2-|\Hess_Qf|_Q^2\bigr)\,d\sigma_Q
 =
 \frac32
 \bigl(|dy^1|_Q^2|dy^2|_Q^2
       -\langle dy^1,dy^2\rangle_Q^2\bigr)
 \Vol_Q(\T^{n-1})>0,
\]
also with a uniform positive lower bound.

For each sufficiently large fixed $L$, choose the
majorant $a$ in the preceding construction independently
of $Q,s$, using the uniform bounds on a compact interval
and exponential decay on the end.
Choose
\[
 b_0(t)=Cs^2e^{-\lambda L}\eta_{0,L}(t),
\]
where $\eta_{0,L}$ is a fixed smooth nonnegative cutoff,
equal to one on $[L,L+1]$ and supported in $(L-1,L+2)$.
A uniform choice of $C$ ensures
$c_n\alpha b_0\ge2R(\bar g)_-$.
Keep the metric cutoff and $\eta_L$ independent of $Q,s$.
The equation
\[
 \int_0^\infty y_\beta(t)\,dt=c(Q,s)^{-1}-1
\]
then determines $\beta$ smoothly in $(Q,s)$, since its
left side is affine in $\beta$ with a fixed positive slope.
The estimates already proved give $\beta\ge0$ for
uniformly small $|s|$ and sufficiently large $L$.
At $s=0$, both $b_0$ and $\beta$ vanish.
Thus the restored metrics $\widehat g_{Q,s,L}$ form
a smooth family, including
$\widehat g_{Q,0,L}=dt^2+Q$, and \eqref{eq:cost}
holds uniformly in $Q$.

Choose $s_*>0$ sufficiently small and fix $L=s_*^{-2}$. The value of $L$ is fixed throughout the family and does
not vary with $s\in[-s_*,0]$. The cubic expansion and \eqref{eq:cost} yield
\[
 \int_{\partial M}
 H_{\widehat g_{Q,-s_*,L}}\,d\sigma_Q
 \ge c_*s_*^3-Cs_*^4
       -Cs_*^2e^{-\lambda/s_*^2}>0
\]
uniformly in $Q$, for some $c_*>0$.
Set $K_{Q,s}=\widehat g_{Q,s,L}$ and $A=L+3$.
This proves the remaining assertions, including
\eqref{eq:family-gain}.
The derivative bounds follow from smooth dependence
on the larger parameter neighborhood and the common
support bound.
\end{proof}

We now localize the periodic models in a boundary chart.
The following proposition quantifies the gain in total
mean curvature and the scalar-curvature error introduced
by localization.
\begin{proposition}\label{prop:localization}
Let $(\Omega,h)$ be a smooth compact connected $n$-manifold
with nonempty boundary, $n\ge3$, and $R(h)\ge0$.
Write $\gamma=h|_{T\partial\Omega}$.
For all sufficiently large integers $N$, there are smooth
metrics $g_N$, uniformly equivalent to $h$, and constants
$a,A,C>0$ independent of $N$ such that
\[
 g_N|_{T\partial\Omega}=\gamma,\qquad
 \F(g_N)\geq aN-C,\qquad
 R(g_N)_-\le CN\mathbf1_{E_N},
\]
where $E_N$ lies in a boundary collar of thickness $A/N$.
The perturbations can be supported in the collar over
any prescribed nonempty open subset of $\partial\Omega$.
\end{proposition}

\begin{proof}
Choose inward Gaussian collar coordinates over a boundary
chart $U\Subset\mathbb R^m$, $m=n-1$, contained in the
prescribed open subset:
\[
 h=dt^2+\gamma_t(x),\qquad
 Q(x)=\bigl((\gamma_0)_{ij}(x)\bigr).
\]
After shrinking $U$, assume that $Q(\overline U)$ lies
in the matrix neighborhood of Lemma~\ref{lem-periodic-family}.
Choose a coordinate cube $D\Subset U$ and choose a smooth $s$ taking values in $[-s_*,0]$ so that
$s=-s_*$ on $D$, with $s=0$ outside a small neighborhood
of $\overline D$;
its choice will be specified below.

We lift the coefficients of $K_{Q,s}$ to
$[0,\infty)\times\mathbb R^m$, where they are
$2\pi$-periodic in each tangential variable, and define
\begin{equation}\label{eq:transplant}
 (g_N)_{ab}(t,x)
 =(K_{Q(x),s(x)})_{ab}(Nt,Nx)+B_{ab}(t,x),
\end{equation}
where index zero denotes the normal coordinate and
\[
 B_{00}=B_{0i}=0,\qquad
 B_{ij}=(\gamma_t)_{ij}-Q_{ij}.
\]

By \eqref{eq:family-end}, $g_N=h$ wherever $s=0$ or
$t\ge A/N$, so it extends smoothly by $h$ outside the chart.
At $t=0$, its induced metric is exactly $\gamma$.
On the layer $0\le t\le A/N$, we have $B=O(N^{-1})$;
hence $g_N$ is positive definite and uniformly equivalent
to $h$ for large $N$.

Differentiating the coefficients in \eqref{eq:transplant},
each derivative in $(Nt,Nx)$ contributes a factor $N$,
whereas derivatives of $Q$ and $s$ are bounded independently
of $N$. Write $z=(\tau,y)$ for the model coordinates.
For each tangential index $i=1,\ldots,m$, the chain rule gives
\begin{equation}
\begin{aligned}
 \partial_i\!\left[K_{Q(x),s(x)}(Nt,Nx)\right]
 &=N(\partial_{y_i}K_{Q(x),s(x)})(Nt,Nx)\\
 &\quad+(D_QK_{Q(x),s(x)})(Nt,Nx)[\partial_iQ(x)]\\
 &\quad+(\partial_sK_{Q(x),s(x)})(Nt,Nx)\,\partial_i s(x)\\
 &=
 N(\partial_{y_i}K_{Q(x),s(x)})(Nt,Nx)+O(1).
\end{aligned}
\end{equation}
Here $D_Q$ denotes differentiation with respect to the
symmetric matrix parameter $Q$, and the brackets indicate
evaluation on the direction $\partial_iQ(x)$. The last two terms are
uniformly bounded because $Q$ and $s$ are fixed smooth
functions and the model family has uniformly bounded
parameter derivatives. Since $Q$ and $s$ are independent
of $t$, we also have
\[
 \partial_t\!\left[K_{Q(x),s(x)}(Nt,Nx)\right]
 =N(\partial_\tau K_{Q(x),s(x)})(Nt,Nx).
\]
Similarly, applying the chain rule once more gives
\begin{equation}\partial_a\partial_b
 \!\left[K_{Q(x),s(x)}(Nt,Nx)\right]=
 N^2(\partial_{z^a}\partial_{z^b}
 K_{Q(x),s(x)})(Nt,Nx)+O(N).
\end{equation}
Here $a,b=0,\ldots,m$, $\partial_0=\partial_t$,
and derivatives in $z$ hold $Q,s$ fixed.
Moreover,
\[
 B=O(N^{-1}),\qquad |DB|+|D^2B|\le C.
\]
All estimates are componentwise and uniform on the layer.
Since $B=O(N^{-1})$ and the model metrics are uniformly
positive definite, the inverse metrics also differ by
$O(N^{-1})$. Multiplying this difference by the
$O(N^2)$ second-derivative and quadratic first-derivative
terms contributes only $O(N)$. Since scalar curvature is linear in second derivatives
and quadratic in first derivatives of the metric, in $U$ we have
\begin{equation}\label{eq:transplantR}
 R(g_N)(t,x)
 =N^2R(K_{Q(x),s(x)})(Nt,Nx)+O(N)\ge-CN.
\end{equation}
Outside $U$, $R(g_N)=R(h)\ge0$, so \eqref{eq:transplantR} holds in $\Omega$.

At the boundary, $B=0$. Since mean curvature is linear
in the first derivatives of the metric, the same argument
gives
\begin{equation}\label{eq:transplantH}
 H_{g_N}(x)
 =N H_{K_{Q(x),s(x)}}(Nx)+O(1).
\end{equation}

We estimate the boundary integral directly using the
positive total mean curvature of the model.
Write $T=\operatorname{supp}s\setminus D$ for the
transition region, whose size will be fixed below.
All volumes $|D|,|T|$ in this argument are measured
with respect to the coordinate measure $dx$.

Consider the cubes
\[
 C_{j,N}=\frac{2\pi}{N}\bigl(j+[0,1]^m\bigr),
 \qquad j\in\mathbb Z^m,
\]
contained in $D$, and choose $x_{j,N}\in C_{j,N}$.
Put $Q_j=Q(x_{j,N})$.
Smooth parameter dependence gives a uniform bound
$|H_{K_{Q,-s_*}}(y)|\le C$.
Together with $|Q(x)-Q_j|\le C/N$ on $C_{j,N}$,
it also yields
\[
 \left|
 H_{K_{Q(x),-s_*}}(Nx)\sqrt{\det Q(x)}
 -
 H_{K_{Q_j,-s_*}}(Nx)\sqrt{\det Q_j}
 \right|\le \frac CN.
\]
Integration over $C_{j,N}$, whose volume is
$(2\pi/N)^m$, gives an error of order $N^{-m-1}$. Smooth dependence on $Q$ and the bounded derivatives
of $Q(x)$ imply
\[
\begin{aligned}
 &\quad \int_{C_{j,N}}
 H_{K_{Q(x),-s_*}}(Nx)\sqrt{\det Q(x)}\,dx\\
 &\quad=
 \sqrt{\det Q(x_{j,N})}
 \int_{C_{j,N}}H_{K_{Q(x_{j,N}),-s_*}}(Nx)\,dx
 +O(N^{-m-1})\\
 &\quad=
 N^{-m}\F(K_{Q(x_{j,N}),-s_*})+O(N^{-m-1})\\
 &\quad\ge bN^{-m}-CN^{-m-1}.
\end{aligned}
\]
where $b$ is from \eqref{eq:family-gain}.
There are $\frac{N^m|D|}{(2\pi)^m}+O(N^{m-1})$
such cubes, and the uncovered part of $D$ has coordinate
volume $O(N^{-1})$.
The uniform bound for the model mean curvatures therefore
gives
\[
 \int_D H_{K_{Q(x),-s_*}}(Nx)\,d\sigma_\gamma
 \ge \frac{b|D|}{(2\pi)^m}-O(N^{-1}).
\]

On the transition region, the same uniform bound gives
\[
 \int_T H_{K_{Q(x),s(x)}}(Nx)\,d\sigma_\gamma
 \ge -C_0|T|,
\]
where $C_0$ is independent of the cutoff derivatives.
Outside $\operatorname{supp}s$, the model is the product
metric and its boundary mean curvature vanishes.
Choose the transition region sufficiently small that
\[
 C_0|T|<\frac{b|D|}{2(2\pi)^m}.
\]
It follows from \eqref{eq:transplantH} that
\begin{equation}\label{eq:averaged-gain}
 \F(g_N)\ge aN-C,\qquad
 a=\frac{b|D|}{2(2\pi)^m}>0.
\end{equation}
The cutoff $s$ is now fixed before $N\to\infty$;
its derivatives only affect the constants in the
preceding curvature estimates.
This completes the proof.

\end{proof}

The following lemma controls the conformal correction
when the negative part of the scalar curvature is
concentrated in a thin boundary layer.

\begin{lemma}\label{lem:repair}
Let $g$ be a fixed smooth metric on a compact connected manifold
$\Omega$ with nonempty boundary, and fix a boundary collar coordinate
$t\ge0$ with $t=0$ on $\partial\Omega$ in a collar neighborhood $U$. Suppose $g_N$ are smooth
metrics, uniformly equivalent to $g$, with common induced boundary
metric $\gamma$. Assume there are constants $A,C>0$ independent of $N$
and measurable sets $E_N$ in the collar such that
\begin{equation}\label{eq:negative-layer}
 R(g_N)_-\le CN\mathbf1_{E_N},\qquad
 E_N\subset\{0\le t\le A/N\}.
\end{equation}
For all sufficiently large $N$, the problem
\begin{equation}\label{thin-layer-conformal}
 (-c_n\Delta_{g_N}+R(g_N))u_N=0,\qquad
 u_N|_{\partial\Omega}=1
\end{equation}
has a unique smooth positive solution. The metrics $\widetilde g_N=u_N^qg_N$
satisfy
\[
 R(\widetilde g_N)=0,\qquad
 \widetilde g_N|_{T\partial\Omega}=\gamma,\qquad
 \F(\widetilde g_N)\ge\F(g_N)-C',
\]
where $C'$ is independent of $N$.
\end{lemma}

\begin{proof}
All integrals and gradient norms in this proof use $g_N$. For smooth
$u$ vanishing on the boundary, we have $
 |u(t,x)|^2\le t\int_0^t|\partial_r u(r,x)|^2\,dr$ near the collar neighborhood of the boundary. Integrating once more in $t$ gives $\int_0^\delta |u(t,x)|^2\,dt\leq \frac{\delta^2}{2}
   \int_0^\delta|\partial_ru(r,x)|^2\,dr$ for small $\delta$, and then uniform
metric equivalence of $g_N$ gives
\begin{equation}\label{eq:layer-poincare}
 \int_{E_N}u^2\le\frac{C_1}{N^2}\int_\Omega|\nabla u|^2.
\end{equation}
For general $u\in H_0^1(\Omega)$, the same inequality holds by a standard approximation argument. In particular,
no derivative bounds on $g_N$ enter \eqref{eq:layer-poincare}.

The Dirichlet quadratic form therefore satisfies
\begin{equation}\label{eq:coercivity}
 \begin{split}
 Q_N(u,u)&\coloneqq c_n\int_\Omega|\nabla u|^2+
                   \int_\Omega R(g_N)u^2\\
 &\ge\left(c_n-\frac{C_2}{N}\right)
                   \int_\Omega|\nabla u|^2.
 \end{split}
\end{equation}
Uniform equivalence also provides a uniform global Dirichlet Poincar\'e
inequality. For sufficiently large $N$, the form is coercive.
Standard elliptic theory gives existence, uniqueness, smoothness and positivity of the solution to \eqref{thin-layer-conformal}.

Set $w_N=(u_N-1)_+$. Testing the equation with $w_N$ gives
\begin{equation}\label{eq:excess-energy}
 Q_N(w_N,w_N)=-\int_\Omega R(g_N)w_N
                 \le\int_{E_N}R(g_N)_-w_N.
\end{equation}
The volume of $E_N$ is $O(N^{-1})$. By \eqref{eq:layer-poincare},
\[
 \int_{E_N}R(g_N)_-w_N
 \le CN\Vol(E_N)^{1/2}\|w_N\|_{L^2(E_N)}
 \le C_3N^{-1/2}\|\nabla w_N\|_{L^2(\Omega)}.
\]
Combining this with \eqref{eq:coercivity}--\eqref{eq:excess-energy}
yields
\[
 \|\nabla w_N\|_{L^2}=O(N^{-1/2}),\qquad
 \int_\Omega R(g_N)_-w_N=O(N^{-1}).
\]
Since $u_N\le1+w_N$, it follows that
\begin{equation}\label{eq:negative-cost}
 \int_\Omega R(g_N)_-u_N
 \le\int_\Omega R(g_N)_-+
      \int_\Omega R(g_N)_-w_N\le C_4.
\end{equation}

Formula \eqref{thin-layer-conformal} gives scalar-flatness. The boundary value
of $u_N$ preserves the induced metric. Integrating the equation and
using \eqref{eq-conformal-RH}, we obtain the exact identity
\begin{equation}\label{eq:repair-cost}
 \F(u_N^qg_N)-\F(g_N)
 =\frac{c_n}{2}\int_{\partial\Omega}\partial_{\nu_{g_N}}u_N
 =\frac12\int_\Omega R(g_N)u_N.
\end{equation}
Positivity of $u_N$ and \eqref{eq:negative-cost} prove the lower bound.
\end{proof}

We now combine Proposition~\ref{prop:localization} and Lemma~\ref{lem:repair} to prove the main theorem.
\begin{proof}[Proof of Theorem~\ref{thm:universal}]
By the boundary extension theorem of Shi--Wang--Wei
\cite[Theorem~1.1]{SWW}, choose a smooth metric $h$ on
$\Omega$ with positive scalar curvature and
$h|_{T\partial\Omega}=\gamma$.
Proposition~\ref{prop:localization} gives metrics $g_N$,
uniformly equivalent to $h$, such that
\[
 g_N|_{T\partial\Omega}=\gamma,\qquad
 \F(g_N)\ge aN-C,\qquad
 R(g_N)_-\le CN\mathbf1_{E_N},
\]
where $a>0$ and $E_N$ lies in a boundary layer of
thickness $A/N$.
Applying Lemma~\ref{lem:repair}, we obtain smooth positive
functions $u_N$, equal to $1$ on $\partial\Omega$, such that
$\widetilde g_N=u_N^qg_N$ satisfies
\[
 R(\widetilde g_N)=0,\qquad
 \widetilde g_N|_{T\partial\Omega}=\gamma,\qquad
 \F(\widetilde g_N)\ge aN-C-C'\longrightarrow+\infty.
\]
This proves the theorem.
\end{proof}

\noindent\textbf{Competing interests.}
The author declares no competing interests relevant to the
content of this article.

\medskip
\noindent\textbf{Data availability.}
Data sharing is not applicable to this article as no datasets
were generated or analysed during the current study.

\medskip
\noindent\textbf{Acknowledgement}
The author is partially supported by the China Postdoctoral Science Foundation (Grant No.~2025T180838).

\medskip
\noindent\textbf{Declaration of AI use}
The author used ChatGPT to assist with revising the manuscript and searching the literature. ChatGPT also assisted in developing the proof of Lemma~\ref{lem:restore} and ~\ref{lem-periodic-family}. The author reviewed the AI-assisted material and takes full responsibility for the content of the paper.

\end{document}